\newcommand*{\memfontfamily}{pazo}
\newcommand*{\memfontpack}{mathpazo}

\documentclass[10pt,%extrafontsizes,%
onecolumn,oneside,a5paper,article%,nofootinbib,nobibnotes,showpacs,amsmath,amssymb,floatfix%,preprintnumbers%tbtags
,frenchb,italian,german,swedish,latin%
,british%
]{memoir}

\usepackage[mathpazo% loads mathpazo - Palatino text & math
,eugreek% greek letters in Euler font
,optima% Zapf's Optima font as sans serif
,itsans% math sans serif in italic
,lmodern% use Latin Modern instead of CM
]{memmana2}

  \hypersetup{
breaklinks=true,%
    pdftitle = {On two recent conjectures in convex geometry},
    pdfauthor = {PierGianLuca  Porta Mana}
  }
\usepackage{mathdots}
\usepackage{breakurl} % uncomment for arxiv

\renewcommand{\multicitedelim}{\addsemicolon\space}
\bibliography{manabibliography}
\renewcommand{\bibfont}{\small}

\DeclareBibliographyCategory{1}
\newcommand{\citep}{\parencites}
\renewcommand{\cite}{\citep}

\setlxvchars
\setxlvchars
\settypeblocksize{*}{\lxvchars}{1.618}
\setulmargins{*}{*}{*}
\setlrmargins{*}{*}{*}
\setheaderspaces{*}{*}{*}
\checkandfixthelayout
\providecommand{\firstname}[1]{#1}
\providecommand{\surname}[1]{#1}
\providecommand{\affiliation}[1]{{\footnotesize#1}}
\providecommand{\epost}[1]{{\footnotesize #1}}
\providecommand{\email}[1]{{\textsf{\href{mailto:#1}{#1}}}}
\providecommand{\pacs}[1]{{\footnotesize\textsc{PACS} numbers: #1}}
\providecommand{\msc}[1]{{\footnotesize\textsc{MSC} numbers: #1}}

\ifdraftdoc
\usepackage%[nodayofweek]
{datetime}
\newdateformat{mydate}{\THEDAY\ \monthname[\THEMONTH] \THEYEAR}
\fi

\newenvironment{acknowledgements}{\chapter*{Acknowledgements}\addcontentsline{toc}{section}{Acknowledgements}}{\par}

\chapterstyle{reparticle}
\copypagestyle{manaart}{plain}
\makeheadrule{manaart}{\headwidth}{0.75\normalrulethickness}
\makeoddhead{manaart}{%
\footnotesize\textit{Porta Mana \amp\ Lewis}}{}{% ***
\footnotesize\textit{On two conjectures in convex geometry}}% ***
\makeoddfoot{manaart}{}{\thepage}{}
\makeatletter
\newcommand\addprintnote{%
\begin{picture}(0,0)%
\put(0,-14){%42
\makebox(0,0){%\rotatebox{0}%
{\tiny%
This document is optimized for on-screen reading and 2-pages-on-1-sheet
printing on A4 or Letter paper}}
}%
\end{picture}%
}
\makeoddfoot{plain}{}{\makebox[0pt]{\thepage}\addprintnote}{%\addNoteLay
}
\makeoddhead{plain}{}{}{\footnotesize
Report no.\ pi-qf-225%***report no.
}

\title{On two recent conjectures in convex geometry%***
}
\author{
{}\\
    \firstname{P.G.L.}\ \surname{Porta\,Mana}\\
\affiliation{Perimeter Institute for Theoretical Physics, Canada}
\\
\epost{\textless\email{lmana AT pitp.ca}\textgreater}
\\[3\jot]
    \firstname{P. G.}\ \surname{Lewis}\\
\affiliation{Quantum Optics and Laser Science Group, Blackett Laboratory,}\\\affiliation{Imperial College, London}
\\
\epost{\textless\email{pgl08 AT imperial.ac.uk}\textgreater}\\{}
}
\predate{\begin{center}\footnotesize}
\date{20 May 2011%\ifdraftdoc\ draft of \mydate\today\fi%
}%***
\postdate{\end{center}}

\usepackage{pifont}

\theoremstyle{plain}

\theoremstyle{remark}

\theoremstyle{definition}

\newcommand{\yC}{\mathcal{C}} %\vdash or \models
\newcommand{\yO}[1]{\mathcal{P}_{#1}}
\newcommand{\yOC}{\yO{\yC}}

\newcommand{\yom}{v}
\newcommand{\yoz}{\yom_0}
\newcommand{\you}{\yom_\text{u}}
\newcommand{\yT}{\varDelta}
\newcommand{\yTT}{\bar{\yT}}
\newcommand{\yOTT}{\yO{\yTT}}

\newcommand{\yf}{F}
\newcommand{\yg}{G}

\newcommand{\yga}{\tfrac{1}{2}}

\newcommand{\yozz}{d_0}
\newcommand{\youu}{d_\text{u}}

\newcommand{\QEM}%{\textnormal{$\Box$}}
{\ding{167}}

\setfloatadjustment{figure}{\footnotesize\centering}
\begin{document}
%\mainmatter
%\bibliographystyle{apsrevmananum} 
%\defaultbibliography{manabibliography}
%\defaultbibliographystyle{apsrevmananum}
\selectlanguage{british}
\hyphenation{ 
Pre-sent Pre-sent-ed Pre-sent-ing Pre-sents Li-ce-o Scien-ti-fi-co Sta-ta-le ca-glia-ri Con-tin-u-um
Quan-tum Be-tween Phe-nom-e-non Mac-ro-scop-ic Mi-cro-scop-ic Sub-space
Sub-spaces Mo-men-tum Mo-men-ta Ther-mo-me-chan-ics Ther-mo-me-chan-i-cal
Meso-scop-ic Elec-tro-mag-net-ic Ve-loc-i-ty Ve-loc-i-ties Gal-i-le-an Gal-i-le-ian
}

%%% Title and abstract %%%
\firmlists*
\maketitle
\abslabeldelim{:\quad}
\setlength{\abstitleskip}{-\absparindent}
\abstractrunin
\begin{abstract}
% ***
  Two conjectures recently proposed by one of the authors are disproved.
  \\[2\jot]
  \msc{52B11,14R99,81P13}\\
  \pacs{02.40.Ft,03.65.Ta,05.90.+m}
\end{abstract}
%\asudedication{*****}

%\setlength{\epigraphwidth}{.7\columnwidth}
%\setlength{\epigraphrule}{0pt}
%\epigraph{`i mean 
%if you dont change in 2 years or even in a day, then whats the point?'}%
%\par\hspace{\stretch{1}}\emph{
%{\textit{L. Pasichnyk}}

\newrefsegment
\selectlanguage{british}

\chapter*{}
\label{cha:into}

The two present authors have found an example of a simplicial refinement
disproving two conjectures recently proposed by one of them (PM)
\citep{portamana2011b}. This example is based on an idea by Rudolph and the
other author (PL) for constructing a `$\psi$-epistemic model' for
finite-dimensional quantum systems \citep{rudolphetal2011}.

The following terminology and symbols are the same as in
ref.~\citep{portamana2011b}.

\bigskip

\begin{figure}[!t]
  \centering
  \includegraphics[width=\columnwidth]{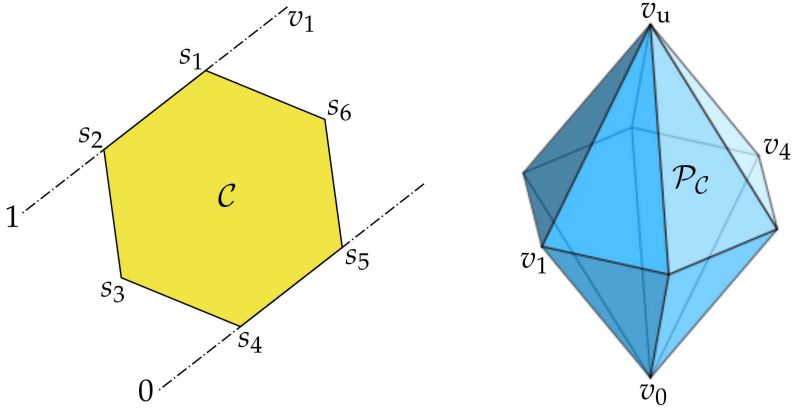}
  \caption{\small The hexagonal convex set $\yC$ and its convex-form space
    $\yOC$. The convex form $v_1$ is represented by two parallel lines on
    $\yC$ and is a point in $\yOC$.}
\label{fig:hexahedron}
\end{figure}

Consider the statistical model given by a hexagon $\yC$ with extreme points
$\set{s_1,\dotsc,s_6}$ and its convex-form space $\yOC$, shown in
fig.~\ref{fig:hexahedron}.
The extreme forms of $\yOC$ are $\set{\yoz,\you,v_1,\dotsc,v_6}$ such that
\pagebreak
\begin{equation}
  \label{eq:extr_v}
{\allowdisplaybreaks  \begin{aligned}
\yoz\inn s_j &=0, \qquad \you\inn s_j =1,
\\     (v_i \inn s_j) &=
    \begin{pmatrix}
  1&1&\yga&0&0&\yga \\
      \yga&1&1&\yga&0&0 \\
      0&\yga&1&1&\yga&0  \\
      0&0&\yga&1&1&\yga  \\
      \yga&0&0&\yga&1&1\\ 
  1&\yga&0&0&\yga&1
    \end{pmatrix},\quad i,j=1,\dotsc,6,
\\
\you&=v_1+v_4=v_2+v_5=v_3+v_6.
  \end{aligned}}
\end{equation}

A simplicial refinement of this statistical model is given by: a
five-dimensional simplex $\yTT$ with extreme points $\set{e_1,\dotsc,e_6}$;
its convex-form space $\yOTT$ (a six-dimensional parallelotope) with $2^6$
extreme points given by the forms
\begin{equation}
  \begin{aligned}
    &\{\yozz, d_1, \dotsc, d_6,\\
 &\qquad d_1+d_2, d_1+d_3, \dotsc, d_5+d_6,\\
      &\qquad\qquad d_1+d_2+d_3, d_1+d_2+d_4, \dotsc, d_4+d_5+d_6,\\
&\qquad\qquad\qquad\dotsc,\\
      &\qquad\qquad\qquad\qquad d_1+\dotsb+d_5,\dotsc,d_2+\dotsb+d_6,\\
%\shoveleft[8em]{
&\qquad\qquad\qquad\qquad\qquad d_1+\dotsb+d_6=\youu\}
%}
  \end{aligned}
\label{eq:etreme_form_6}
\end{equation}
such that
\begin{equation}
\yozz \inn e_j = 0,\qquad
d_i \inn e_j = \delt_{ij};
\label{eq:extr10}
\end{equation}
the partial map $\yf\colon \yTT \hookrightarrow \yC$ defined by
\begin{equation}\label{eq:F_hex}
  \begin{gathered}
    \yf\bigl[\tfrac{1}{2}(e_1+e_2)] = s_1,
\quad   \yf\bigl[\tfrac{1}{2}(e_2+e_3)\bigr] = s_2,
\quad    \dotsc,
\quad    \yf\bigl[\tfrac{1}{2}(e_6+e_1)\bigr] = s_6,
  \end{gathered}
\end{equation}
\begin{figure}[!p]
  \centering
  \includegraphics[width=\columnwidth]{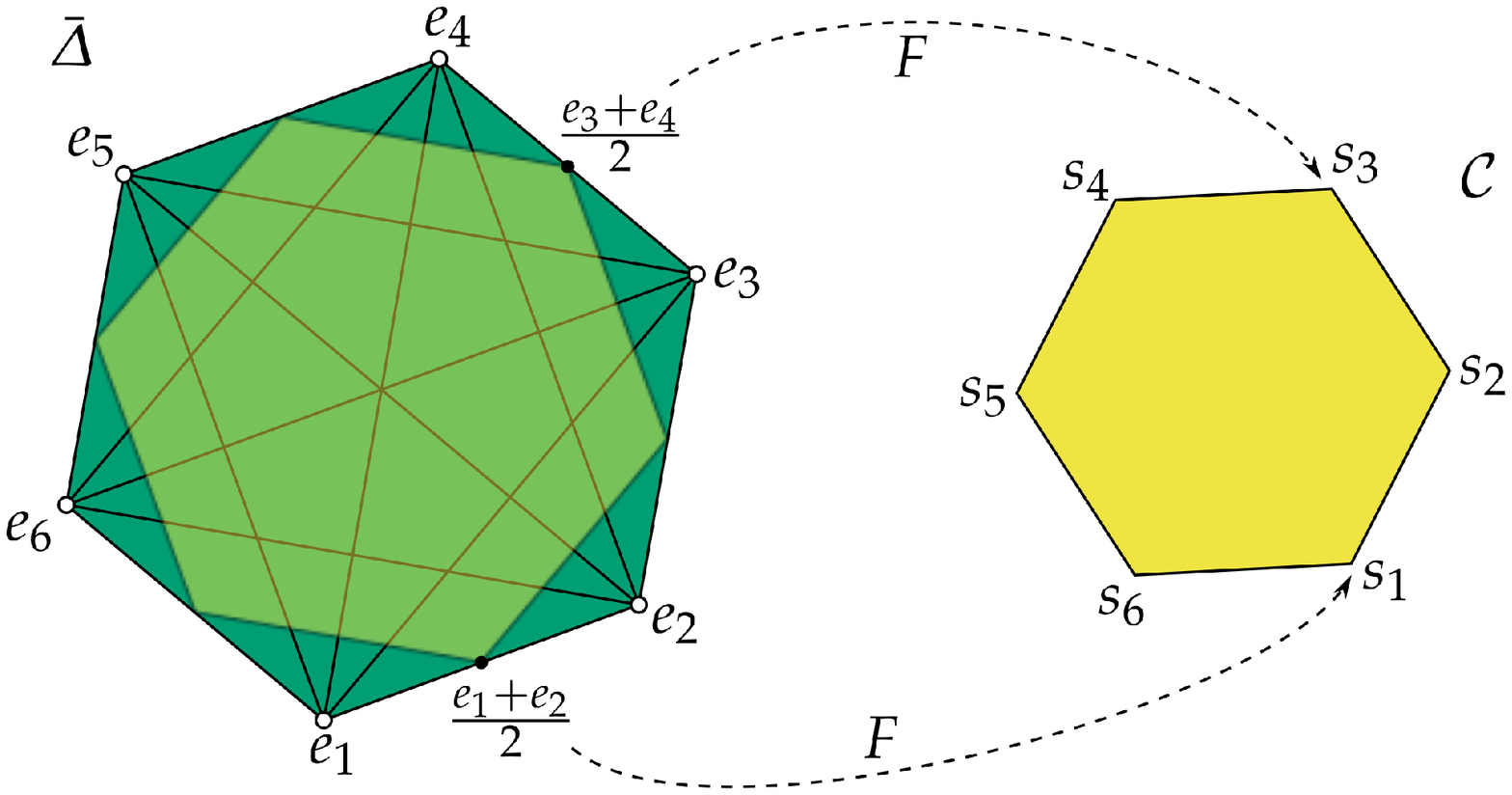}
  \caption{\small Illustration of the partial map $\yf\colon \yTT
    \hookrightarrow \yC$. The domain of definition of this map is a
    four-dimensional polytope with six extreme points, whose projection is
    shown in lighter green.}
\label{fig:hexaref}
\makebox{}\par\vfill\makebox{}
  \includegraphics[width=\columnwidth]{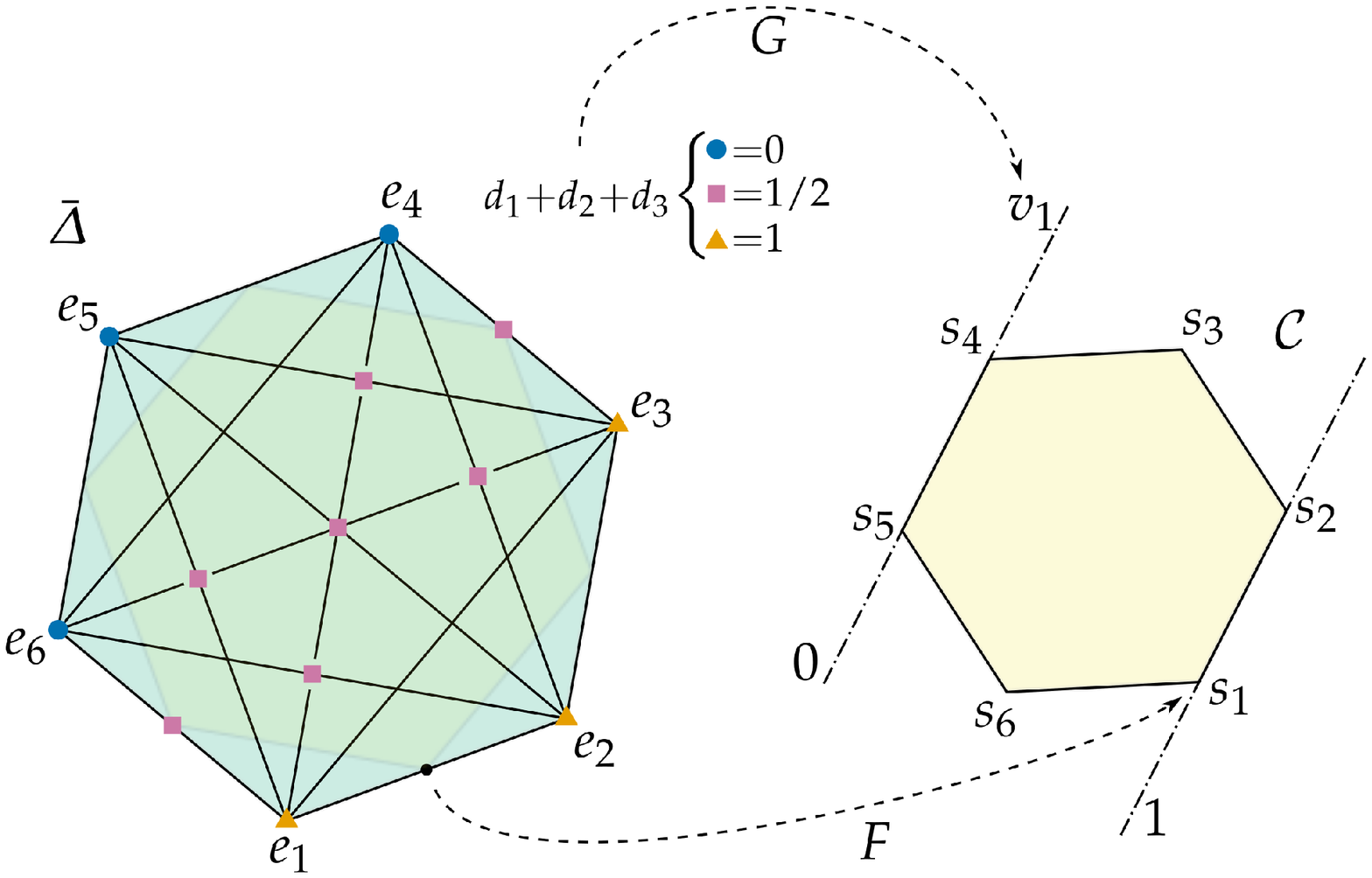}
  \caption{\small Illustration of the partial map $\yg\colon \yOTT
    \hookrightarrow \yOC$. The convex form $d_1+d_2+d_3$ on $\yTT$ is
    defined by its unit value on the extreme points $\set{e_1,e_2,e_3}$
    (orange triangles) and
    zero value on $\set{e_4,e_5,e_6}$ (blue circles); some points (purplish
    squares) on which
    it has value $1/2$ are also shown. This form is mapped by $\yg$ into
    the form $v_1$ of $\yC$.}
\label{fig:hexaform}
\end{figure}
illustrated in fig.~\ref{fig:hexaref}; and finally the partial map
$\yg\colon \yOTT \hookrightarrow \yOC$ defined by
\begin{equation}\label{eq:G_hex}
  \begin{gathered}
\yg(\yozz)=\yoz,\quad
    \yg(d_1+d_2+d_3) = v_1,
\quad    \yg(d_2+d_3+d_4) = v_2,
\quad\dotsc,\\
   \yg(d_5+d_6+d_1) = v_5,
\quad   \yg(d_6+d_1+d_2) = v_6
  \end{gathered}
\end{equation}
(note that $\yg(\youu) = \you$), illustrated in fig.~\ref{fig:hexaform}.

\bigskip

The important features of this refinement are these:
\begin{enumerate}[a.]
\item each extreme point of $\yC$ corresponds to a non-extreme point
  of $\yTT$, and to that alone;
\item each of the eight extreme points of $\yOC$ corresponds to
  an extreme point of $\yOTT$, and to that alone;
\item the map $\yf$ is partial, \ie\ it is not defined on some points of
  $\yTT$, not even its six extreme ones $\set{e_i}$;
  \item\label{item:faces}the minimal faces of $\yTT$ containing the
  counter-images $\yf^{-1}(s_j)$ of the extreme points of $\yC$ have a
  non-empty intersection; for example, the minimal faces containing
  $\yf^{-1}(s_1)$ and $\yf^{-1}(s_2)$ are respectively the edges $e_1e_2$
  and $e_2e_3$, whose intersection is the extreme point $e_2$;
  \item\label{item:proje}the convex subset of $\yTT$ on which $\yf$ is
  defined is a four-\bd dimensional non-simplicial polytope with six
  extreme points; on it, $\yf$ acts like a parallel projection onto $\yC$.
\end{enumerate}
Notably, features~\ref{item:faces}.\ and~\ref{item:proje}.\ respectively
disprove Conjecture~3 and~1 of ref. \citep{portamana2011b}. Conjecture~2 of
the same work, however, still stands.

It is also worth noting that the construction on which this refinement is
based does not work for a square (as shown in the counter-example of
\citep{portamana2011b}) nor, apparently, with a heptagon or an octagon.
Many questions still remain. What can we say about the simplicial
refinements of polyhedra and their convex-form spaces? What general
theorems can we find about simplicial refinements in general?

\begin{acknowledgements}
 Many thanks to the developers and maintainers of \LaTeX, Emacs, AUC\TeX,
  MiK\TeX, arXiv, Inkscape.
  Research at the Perimeter Institute is supported by the Government of
  Canada through Industry Canada and by the Province of Ontario through the
  Ministry of Research and Innovation.
\end{acknowledgements}

%\appendixpage
%\appendix

%%%%%%%%%%%%%%% BIB %%%%%%%%%%%%%%%
\pagebreak[1]
\defbibnote{prenote}{%
\small% Note:
%Of two years separated by a virgule, the first is that of original
%publication or composition.
%\texttt{arXiv} eprints available at \url{http://arxiv.org/}.
%\texttt{mp_arc} eprints available at \url{http://www.ma.utexas.edu/mp_arc/}.
%\\
%\texttt{philsci} eprints at \url{http://philsci-archive.pitt.edu/}.
}

\newcommand{\citein}[2][]{\textnormal{\textcite[#1]{#2}}\addtocategory{1}{#2}}
\newcommand{\citebi}[2][]{\textcite[#1]{#2}\addtocategory{1}{#2}}
%\printbibliography[segment=1,prenote=prenote]
\defbibfilter{1}{\category{1} \or \segment{1}}
%\defbibheading{bibliography}{}
\printbibliography[filter=1,prenote=prenote]
%\printbibliography[segment=1,category=extra,prenote=prenote]

%\newrefsegment

\end{document}